\documentclass[11pt]{article}
\usepackage{times}
\usepackage{mathfont}

\def\exp{\mathop{\rm exp}}
\def\min{\mathop{\rm min}}

\DeclareSymbolFont{AMSb}{U}{msb}{m}{n}
\DeclareSymbolFontAlphabet{\Bbb}{AMSb}
\def\R{\ensuremath{\Bbb R}}

\def\dg{\ensuremath{^\circ}}

\input epsf
\def\efig#1#2{\hbox{\epsfxsize=#1\epsfbox{#2}}}

\newtheorem{theorem}{Theorem}
\newtheorem{lemma}{Lemma}

\newtheorem{defn}{Definition}

\DeclareSymbolFont{lasy}{U}{lasy}{m}{n}
\let\Box\undefined
\DeclareMathSymbol\Box{0}{lasy}{"32}
\newcommand{\qed}{\hfill$\Box$}
\newenvironment{proof}{\noindent{\bf Proof:}}{\qed\medskip}

\makeatletter

\long\def\@makecaption#1#2{
   \vskip 10pt 
   \setbox\@tempboxa\hbox{{\small #1. #2}}
   \ifdim \wd\@tempboxa >\hsize   
       {\small #1. #2}\par        
     \else                        
       \hbox to\hsize{\hfil\box\@tempboxa\hfil}  
   \fi}

\def\@begintheorem#1#2{\it\trivlist
  \item[\hskip\labelsep{\bf #1\ #2.\ }]}
\def\@opargbegintheorem#1#2#3{\it\trivlist
  \item[\hskip\labelsep{\bf #1\ #2\ {\rm(#3)}.}]}

\makeatother

\begin{document}
\bibliographystyle{abbrv}

\title{Diameter and Treewidth in Minor-Closed Graph Families}

\author{David Eppstein\thanks{Department of Information and Computer
Science, University of California, Irvine, CA 92697-3425,
eppstein@ics.uci.edu, http://www.ics.uci.edu/$\sim$eppstein/.
Supported in part by NSF grant CCR-9258355 and by matching funds from
Xerox Corp.}}

\date{}
\maketitle   

\begin{abstract}
It is known that any planar graph with diameter $D$ has treewidth
$O(D)$, and this fact has been used as the basis for several planar
graph algorithms. We investigate the extent to which similar relations
hold in other graph families. We show that treewidth is bounded by a
function of the diameter in a minor-closed family, if and only if some
apex graph does not belong to the family.  In particular, the $O(D)$
bound above can be extended to bounded-genus graphs.
As a consequence, we extend several approximation algorithms and exact
subgraph isomorphism algorithms from planar graphs to other graph
families.
\end{abstract}

\section{Introduction}

Baker \cite{Bak-JACM-94} implicitly based several planar graph
approximation algorithms on the following result, which can be found
more explicitly in \cite{Bod-EATCS-88}:

\begin{defn}
A {\em tree decomposition} of a graph $G$ is a representation of $G$
as a subgraph of a chordal graph $G'$.  The {\em width} of the tree
decomposition is one less than the size of the largest clique in $G'$.
The {\em treewidth} of $G$ is the minimum width of any
tree decomposition of $G$.
\end{defn}

\begin{lemma}\label{diam-width}
Let $D$ denote the diameter of a planar graph $G$.  Then a
tree decomposition of $G$ with width $O(D)$ can be found in time $O(Dn)$.
\end{lemma}

The lemma can be proven by defining a chordal graph having cliques for
certain three-leaf subtrees in a breadth first search tree of $G$.
Such a subtree has at most $3D-2=O(D)$ vertices.

Baker used this method to find approximation schemes for the maximum
independent set and many other covering and packing problems in planar
graphs, improving previous results on planar graph approximation
algorithms based on separator decomposition
\cite{LipTar-SJC-80,ChiNisSai-JIP-81}. Baker's basic idea was to remove
the vertices in every $k$th level of the breadth first search tree of an
arbitrary planar graph $G$; there are
$k$ ways of choosing which set of levels to remove, at least one of which
only decreases the size of the maximum independent set by a factor of
$(k-1)/k$. Then, each remaining set of contiguous levels forms a graph
with treewidth $O(k)$ (it is a subgraph of the graph with diameter $k$
formed by removing vertices in outer levels and contracting edges in
inner levels), and the maximum independent set in each such
component can be found by standard dynamic programming techniques
\cite{BerLawWon-Algs-87, TakNisSai-JACM-82}.

Other workers have develeped
parallel variants of these approximation schemes
\cite{Che-ISAAC-95,CheHe-WG-95,DiaSerTor-AI-96}, applied Baker's method to
exact subgraph isomorphism, connectivity, and shortest path algorithms
\cite{Epp-JGAA-97}, extended similar ideas to approximation algorithms in
other classes of graphs
\cite{Che-ICALP-96,ThiBod-IPL-97} or graphs equipped with a geometric
embedding
\cite{HunMarRad-ESA-94}, and definited structural
complexity classes based on these methods
\cite{KhaMot-STOC-96}.

These results naturally raise the question, how much further can these
algorithms be extended? To what other graph families
do these techniques apply?
Since the argument above about
contiguous levels of the breadth first search tree being contained in a
low-diameter graph is implicitly based on the concept of {\em graph
minors}, we restrict our attention to {\em minor-closed} families; that
is, graph families closed under the operations of edge deletion and edge
contraction. 
Minor-closed families have been studied extensively by Robertson,
Seymour, and others, and include such familiar graph families as the
planar graphs, outerplanar graphs, graphs of bounded genus, graphs of
bounded treewidth, and graphs embeddable in $\R^3$ without any linked or
knotted cycles.

\begin{defn}
Define a family $\cal F$ of graphs to have the {\em
dia\-me\-ter-tree\-width property} if
there is some function $f(D)$ such that every graph in $\cal F$ with
diameter at most $D$ has treewidth $f(D)$.
\end{defn}

\begin{figure}
$$\efig{4in}{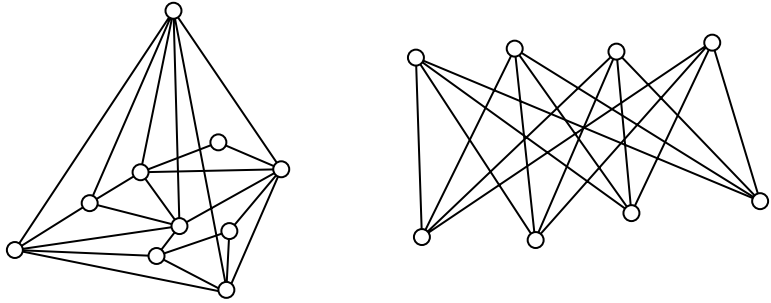}$$
\caption{The graph on the left is an apex graph; the topmost vertex is
one of the possible choices for its apex. The graph on the right
is not an apex graph.}
\label{apex}
\end{figure}

Lemma~\ref{diam-width} can be rephrased as showing that the planar
graphs have the diameter-treewidth property with $f(D)=O(D)$.  In this
paper we exactly characterize the minor-closed families of graphs having
the diameter-treewidth property, in a manner similar to Robertson and
Seymour's characterization of the minor-closed families with bounded
treewidth as being those families that do not include all planar
graphs~\cite{RobSey-GM5}.

\begin{defn}
An {\em apex graph} is a graph $G$ such that for some vertex~$v$
(the {\em apex}), $G-v$ is planar (Figure~\ref{apex}).
\end{defn}

Apex graphs have also been known as nearly-planar graphs,
and have been introduced to study linkless and knotless 3-dimensional
embeddings of graphs~\cite{RobSeyTho-GST-91,Wel-GST-91}.

The significance of apex graphs for us is that they provide examples of graphs
without the diameter-treewidth property: let $G$ be an $n\times n$
planar grid, and let $G'$ be the apex graph formed by connecting some
vertex~$v$ to all vertices of $G$; then $G'$ has treewidth $n+1$
and diameter~2.  Therefore, the family of apex graphs does not have the
diameter-treewidth property, nor does any other family containing all apex
graphs.  Our main result is a converse to this: any minor-closed family
$\cal F$ has the diameter-treewidth property, if and only if $\cal F$
does not contain all apex graphs.

\section{Walls}

Recall that the Euclidean plane can be exactly covered by translates of
a regular hexagon, with three hexagons meeting at a vertex.

\begin{defn}
We say that a set of hexagons is {\em connected}
if its union is a connected subset of the Euclidean plane.
If $h_1$ and $h_2$ are two hexagons from a
tiling of the plane by infinitely many regular hexagons,
define the {\em distance} between the two hexagons to be the smallest
integer $d$ for which there exists a connected subset of the infinite
tiling, containing both $h_1$ and $h_2$, with cardinality $d+1$.
\end{defn}

Thus, any hexagon is at distance zero from itself; two hexagons meeting
edge-to-edge are at distance one, and in general if $h_1\neq h_2$ and
$h_1$ is at distance $d$ from $h_2$, then $h_1$ meets edge-to-edge with
some other hexagon at distance $d-1$ from $h_2$.

\begin{figure}[t]
$$\efig{5in}{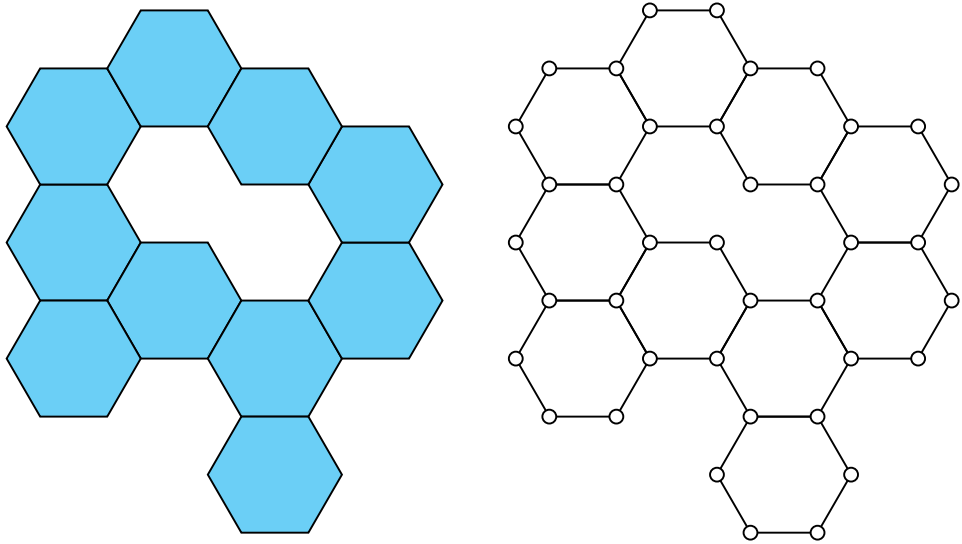}$$
\caption{A set of hexagons and its graph.}
\label{hexgraph}
\end{figure}

\begin{defn}
Let $S$ be a finite connected subset of the hexagons in a
tiling of the Euclidean plane by regular hexagons.
Then we define the {\em
graph of $S$} to be formed by creating a vertex at each point of the
plane covered by the corner of at least one tile of $S$, and creating an
edge along each line segment of the plane forming one of the six edges
of at least one hexagon in $S$.
\end{defn}

Observe that the graph of $S$ is planar and each of
its vertices has degree at most three.
Figure~\ref{hexgraph} shows an example of a set of hexagonal tiles and
its graph.

\begin{figure}[t]
$$\efig{5in}{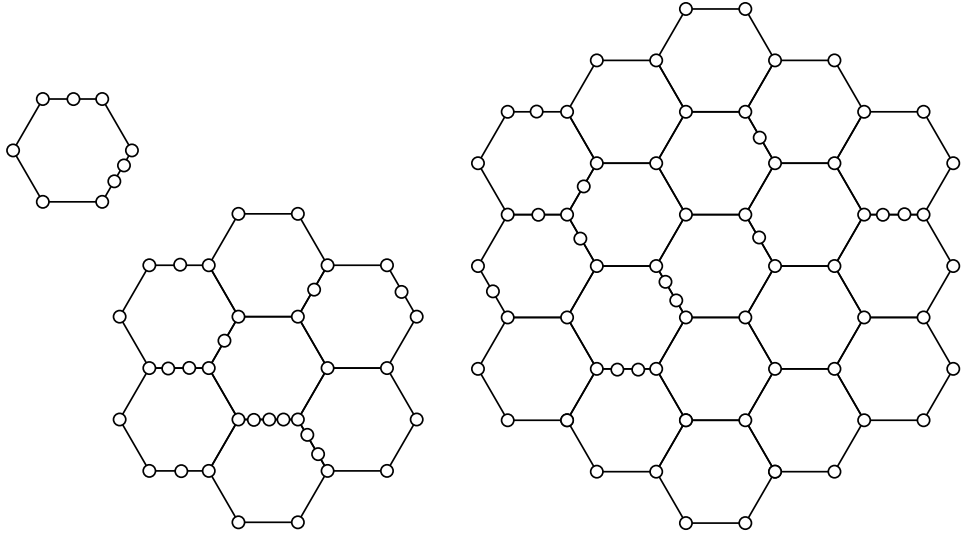}$$
\caption{Walls of size one, two, and three.}
\label{wallfig}
\end{figure}

\begin{defn}
A {\em subdivision} of a graph $G$ is a graph $G'$
formed by replacing some or all edges of $G$ by paths of two or more
edges. A {\em wall of size $s$} is a subdivision of the graph of $S$,
where $S$ is the set of all hexagons within distance
$s-1$ from some given {\em central tile} in a tiling of the Euclidean
plane by regular hexagons.
\end{defn}

Note that since the definition of a wall depends only on the
combinatorial structure of $S$, it is independent of the particular
tiling or central tile chosen in the definition. Examples of
walls are shown in Figure~\ref{wallfig}. Walls are very similar to
(subdivisions of) grid graphs but have a slight advantage of having
degree three.  Thus we can hope to find them as subgraphs rather than as
minors in other graphs.

\begin{lemma}[Robertson and Seymour~\cite{RobSey-GM5}]\label{wall}
For any $s$ there is a number $w=W(s)$ such that any graph of
treewidth $w$ or larger contains as a subgraph a wall of size $s$.
\end{lemma}

In a recent improvement to this lemma, Robertson, Seymour, and Thomas
\cite{RobSeyTho-JCTB-94} showed that if $H$ is a planar graph, the family
of graphs with no
$H$-minor has treewidth at most $20^{2(2|V(H)|+4|E(H)|)^5}$.
Since a wall of size $s$ is a planar graph with $O(s^2)$ edges and
vertices, this implies that $W(s)\le \exp(O(s^{10}))$.

\begin{lemma}[Robertson and Seymour~\cite{RobSey-GM5}]\label{wall-minor}
For any planar graph $G$ there is some $s=s(G)$
such that any wall of size $s$ has $G$ as a minor.
\end{lemma}

We will subsequently need to identify certain components of walls.
To do this we need to use not just the graph-theoretic structure of a
wall but its geometric structure as a subdivision of the graph of a set
of hexagons. (This geometric structure is essentially unique for large
walls, but not for walls of size two, and in any case we will not prove
uniqueness here.)

\begin{defn}
An {\em embedding} of a wall $G$ is the identification of $G$
as a subdivision of a graph of a set of hexagons meeting the
requirements for the definition of a wall. A {\em
$t$-inner} vertex of an embedded wall is a vertex incident to a hexagon
within distance $t-1$ of the wall's central hexagon (so all vertices in a
wall of size $s$ are $s$-inner). An {\em outer} vertex or edge of an
embedded wall is a vertex or edge incident to the boundary of the union of
the set hexagons forming the embedding.
\end{defn}

\section{Routing Across Walls}

\begin{figure}[t]
$$\efig{2in}{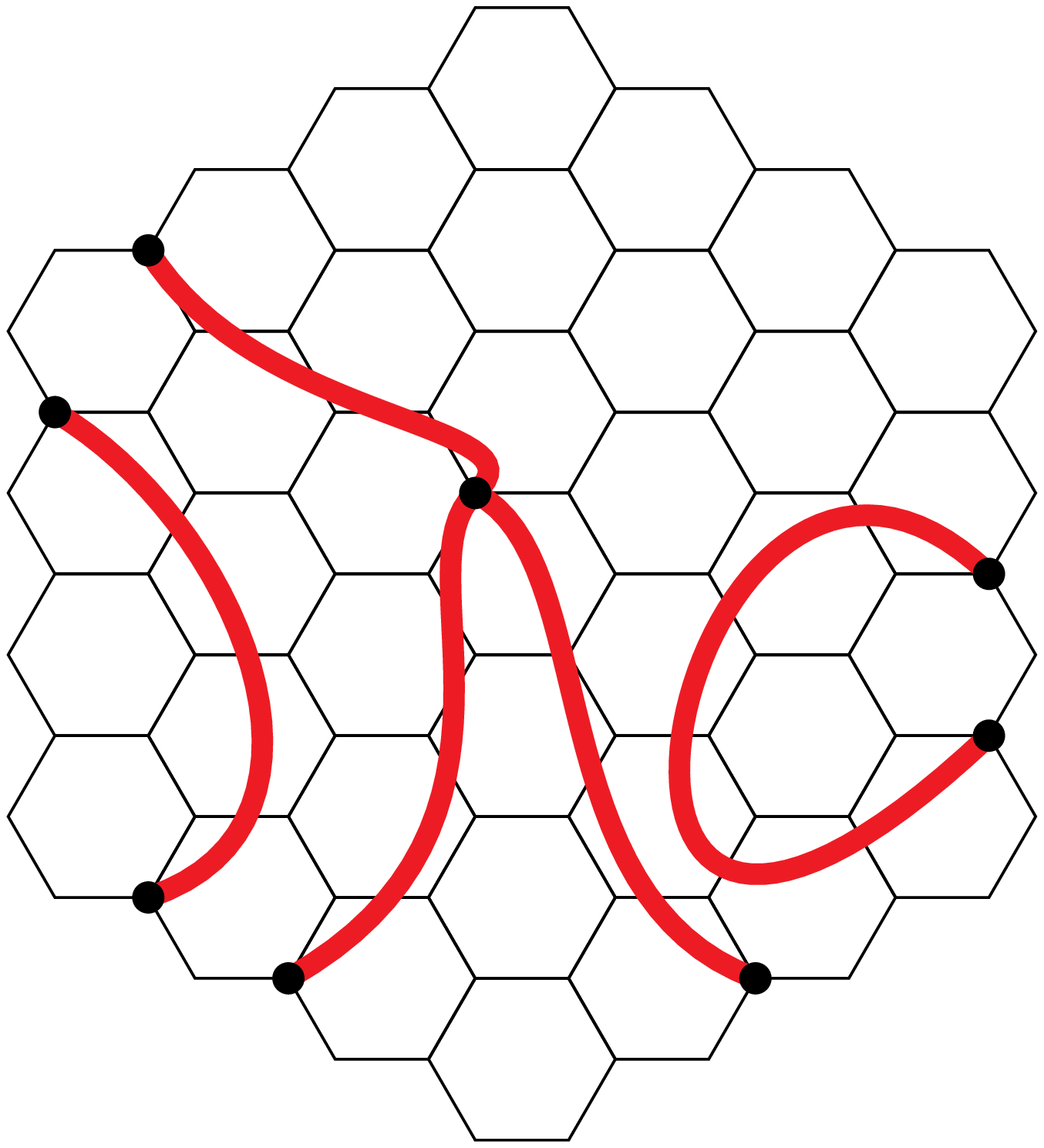}\qquad
\efig{2in}{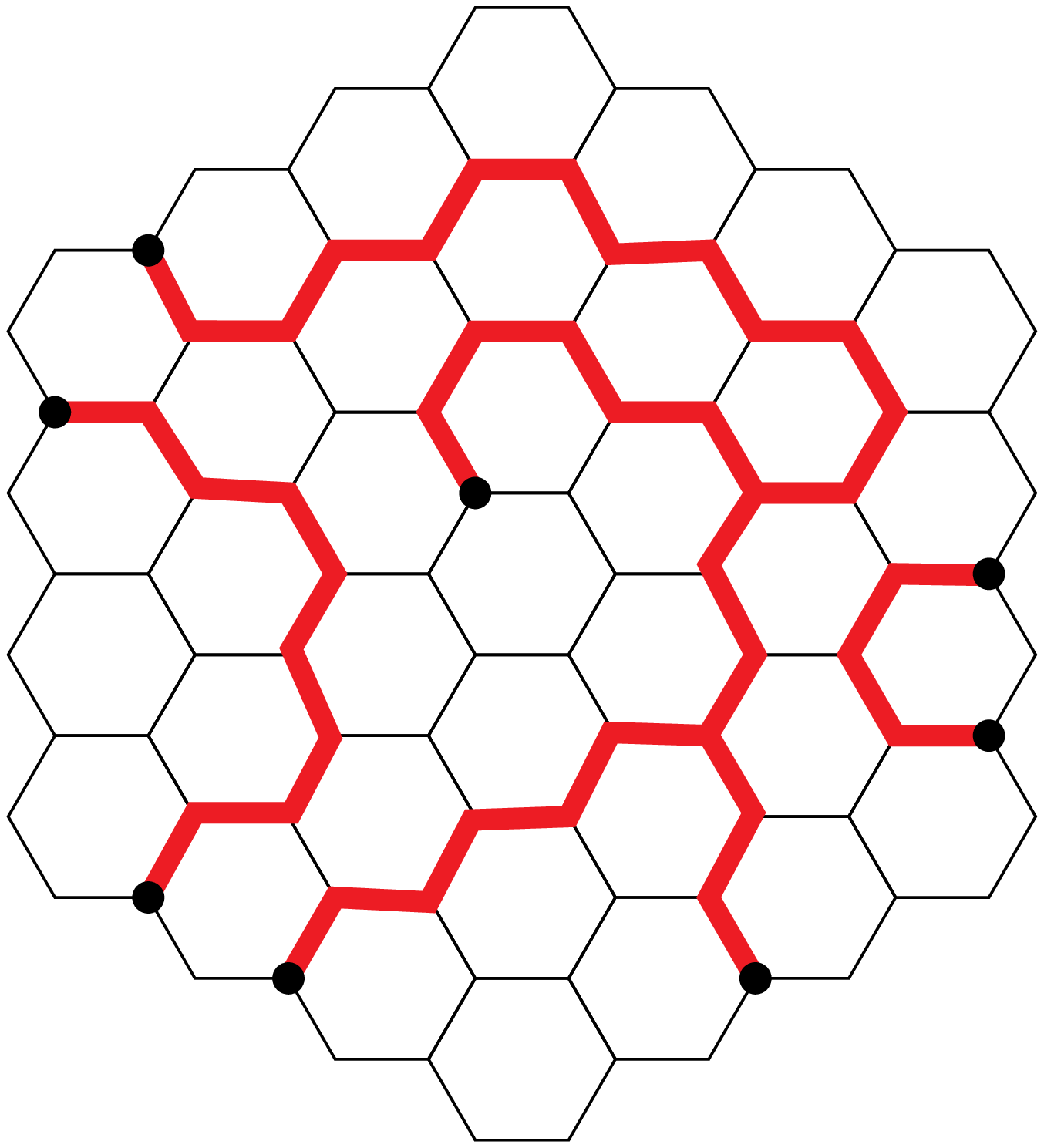}$$
\caption{A $(2,2)$ routing problem and one of its solutions.}
\label{routing}
\end{figure}

\begin{figure}[t]
$$\efig{2in}{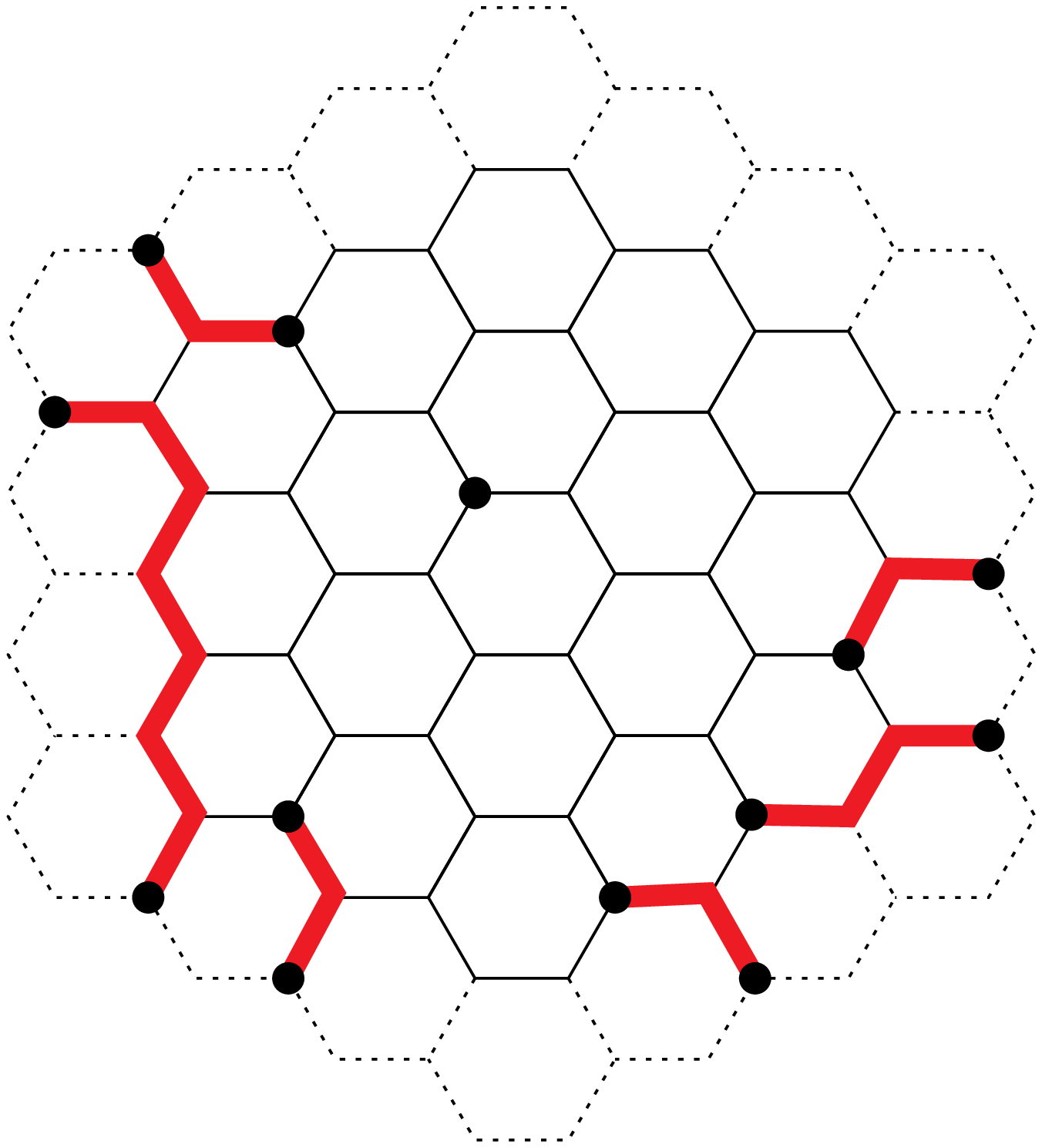}\qquad
\efig{2in}{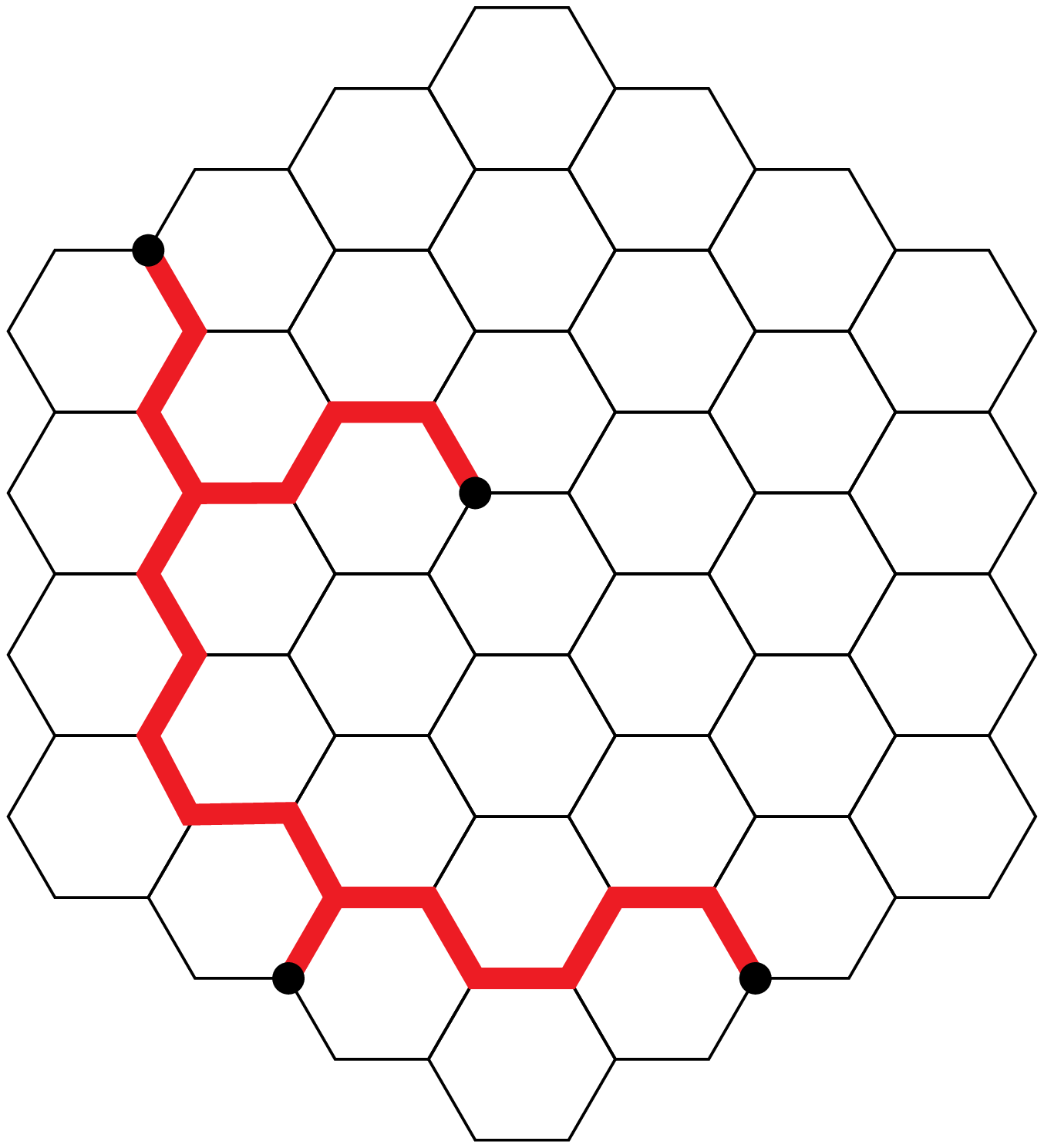}$$
\caption{Cases for solution of routing problems: (a) reduction from an
$(s,t)$ problem to a $(s,t-1)$ problem; (b) solution of $(s,0)$ problem by
routing around boundary of inner wall.}
\label{routed}
\end{figure}

\begin{defn}
An {\em $(s,t)$ routing problem} consists of an embedded wall $G$ of size
$s+t$, an $(s-1)$-inner vertex $v$, and a set $S$ of pairs of {\em
terminals} (certain vertices of the wall), satisfying the following
conditions:
\begin{enumerate}
\item Each terminal is either $v$ or a degree-three outer vertex of the
wall.
\item Each outer vertex occurs at most once as a terminal in $S$;
$v$ occurs at most three times as a terminal.
\item The graph formed by the pairs of terminals in $S$
has a planar embedding as a set of non-crossing curves within the
interior of $U$, where $U$ denotes the union of the hexagons forming the
embedding of the wall.
\item At most $t$ pairs of terminals do not involve~$v$.
\end{enumerate}
\end{defn}

\begin{defn}
A {\em solution} to an $(s,t)$ routing problem
consists of a vertex $v'$ and
a set of $|S|+1$ edge-disjoint paths in $G$,
satisfying the following conditions:
\begin{enumerate}
\item Each pair in $S$ must correspond to one of the paths of the
solution.
\item Each outer terminal of a pair in $S$ must be an endpoint of the
corresponding path.
\item Each pair in $S$ involving vertex $v$ must correspond to a path
having $v'$ as one of its endpoints.
\item The remaining path in the set, not corresponding to a pair in $S$,
must have as its two endpoints $v$ and one of the vertices on a path
involving $v'$.
\item All paths are disjoint from the outer edges of the wall.
\end{enumerate}
\end{defn}

A $(2,2)$ routing problem (with five paths, three involving the inner
vertex) and its six-path solution is depicted in Figure~\ref{routing}.

\begin{defn}
A pair $(x,y)$ of terminals in an $(s,t)$ routing problem is {\em
splittable} if 
the curve corresponding to $(x,y)$ in the planar
embedding of $S$ partitions
$U$ into two regions $A$ and $B$ such that all terminals are incident to
$A$ and only terminals $x$ and $y$ are incident to $B$.
\end{defn}
In Figure~\ref{routing}, both pairs of outer terminals are splittable.

\begin{lemma}
If an $(s,t)$ routing problem includes a pair of outer terminals, it
includes a splittable pair.
\end{lemma}

\begin{proof}
The planar embedding of the pairs of outer terminals in $S$, together with
the boundary of the wall, forms an outerplanar graph (a planar graph in
which all vertices are incident to the outer face).
Because the weak dual of an outerplanar graph (the graph formed from the
planar dual by removing the vertex corresponding to the outer face) is a
tree, it has at least two leaves.
Each leaf of this tree corresponds to a region of $U$ bounded by a
curve in the planar embedding of $S$ and not containing any other outer
terminals of $S$.  At most one leaf contains the inner terminal $v$,
so there is at least one leaf not containing any terminals.
\end{proof}

\begin{lemma}
Every $(s,t)$ routing problem has a solution.
\end{lemma}

\begin{proof}
We use induction on $t$.  If there are fewer than $t$ pairs of terminals
involving $v$, the given problem is also an $(s+1,t-1)$ routing problem
and the result follows from induction.

If $t>0$, let $(x,y)$ be a splittable pair.
Then we can extend an edge from each terminal of $S$ to an
$(s+t-1)$-inner vertex, on the boundary of a wall of size $(s+t-1)$
within the original wall.  We connect $x$ and $y$ by a path around the
boundary of this smaller wall.

Next we connect each other outer terminal to a
degree-three outer vertex of the smaller wall, one terminal at a time,
starting from the terminal immediately counterclockwise of the pair
$(x,y)$, and continuing counterclockwise from there.
For each terminal $t$, we first attempt to extend a path clockwise around
the inner wall's boundary to the next degree-three vertex.
There are three possible situations that can arise in this extension:
\begin{enumerate}
\item We reach an unused degree-three vertex. This vertex will become
the terminal of a smaller problem in the inner wall.
\item We reach a vertex that is part of a path extended from the other
endpoint $u$ of a pair $(t,u)$ in $S$. In this case we have found a path
connecting $(t,u)$ and will not continue using this pair in the smaller
problem we form.
\item We reach a vertex that is part of a path extended from another
terminal $u$, and both $(t,v)$ and $(u,v)$ are pairs in $S$.
In this case we will form a smaller problem in which these two pairs
have been replaced by a single pair $(w,v)$ where $w$ is the
degree-three vertex reached from both $t$ and~$u$.
\item The degree-three vertex we reach is already part of a path but can
not be connected to $t$. In this final case we instead extend a path
counterclockwise from $t$ to the next degree-three vertex.
\end{enumerate}
Note that the first time the counterclockwise extension of case~4
happens can only be at one of the six points where two degree-two
outer vertices of the wall are adjacent. Case~4 may then continue to
happen as long as each successive degree-three
vertex on the boundary of the wall is a terminal that can not be
connected to the previous terminal.  But, by planarity, this can only
happen if no two terminals in this sequence form pairs with each other or
with $v$, for if they did we would have one of cases 2 or~3 instead. 
Therefore there are at most
$t$ terminals in such a sequence, and we will escape from this
counterclockwise case before we reach the next pair of two adjacent outer
degree-two vertices of the wall. As a consequence, this case always
succeeds in extending the path to an unused degree-three vertex.

The result of this path extension process is an $(s,t-1)$ routing problem
on the smaller wall (Figure~\ref{routed}(a)). By induction, this smaller
problem has a solution which can be combined with the path extensions to
solve the original $(s,t)$ routing problem.

Finally, if $t=0$, we have at most three outer terminals on the boundary
of a wall of size $s$ and one non-boundary vertex $v$.
Again, we extend an edge from each terminal to a vertex on the boundary
of a smaller wall of size $(s+t-1)$.  We connect these three vertices by
paths. If $v$ is not already on one of these paths
we add a path connecting it to the solution (Figure~\ref{routed}(b)).
\end{proof}

\section{Macrocells}

\begin{figure}[t]
$$\efig{3in}{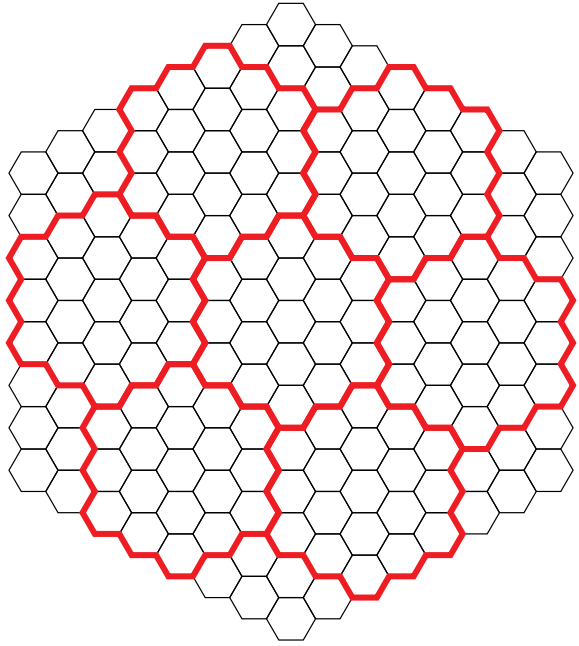}$$
\caption{Subdivision of a large wall into many smaller walls.}
\label{regions}
\end{figure}

\begin{figure}[t]
$$\efig{5in}{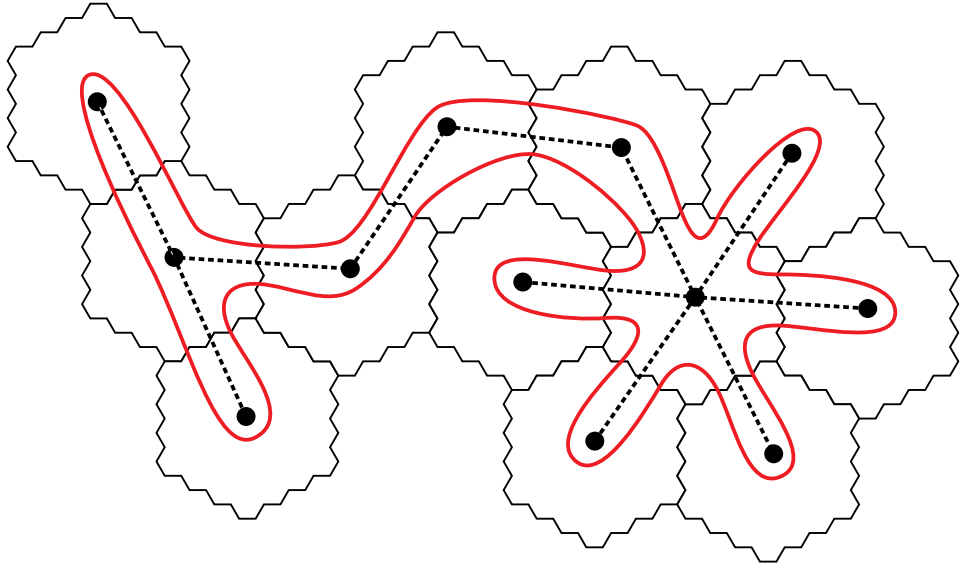}$$
\caption{Curve formed by doubling the spanning tree of a set of
macrocells.}
\label{wallcurve}
\end{figure}

\begin{figure}[t]
$$\efig{5in}{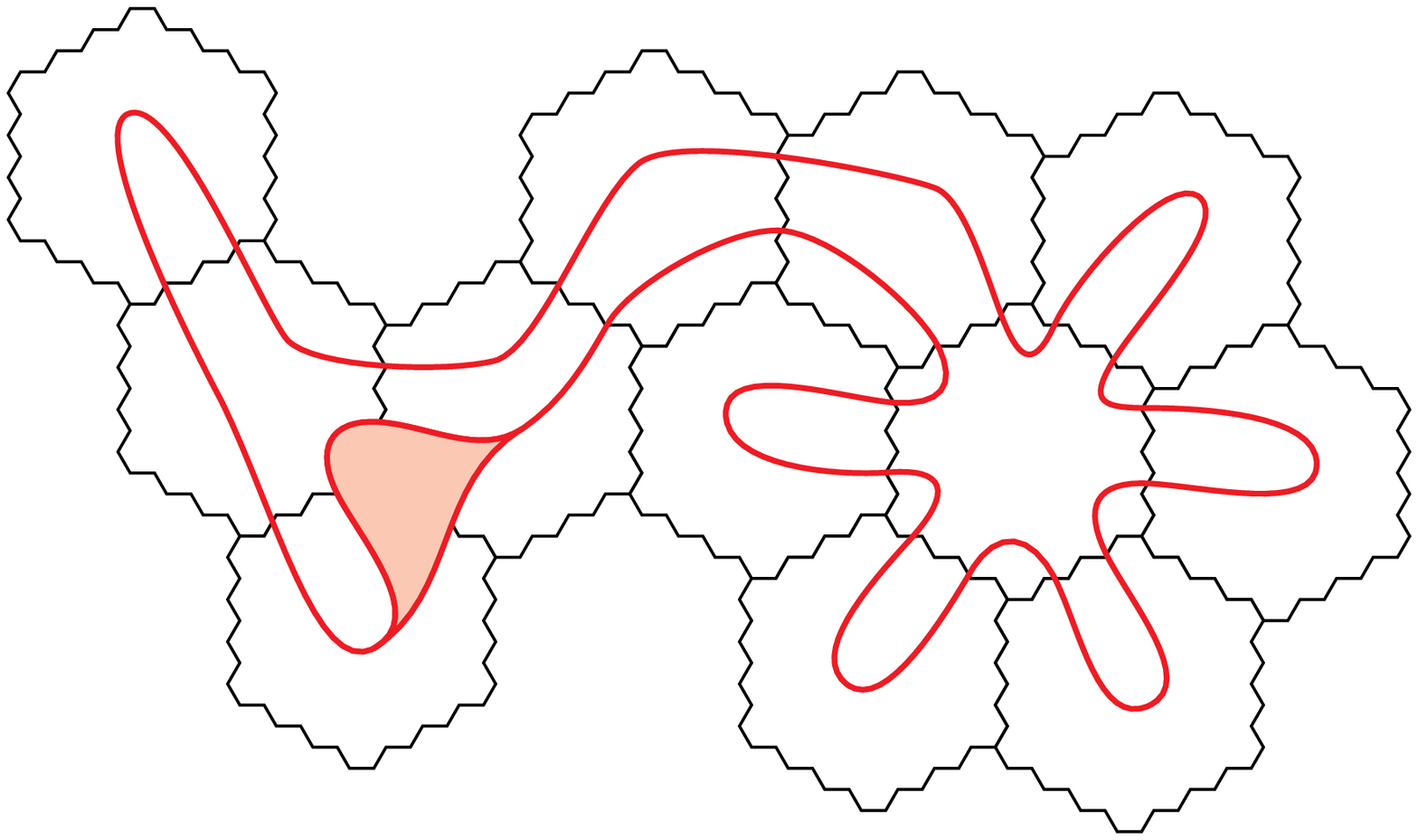}$$
\caption{Curve simplification by removal of a $60\dg$ angle.}
\label{curvesimp}
\end{figure}

The strategy for our proof that graphs without the diameter-treewidth
property contain all apex graphs as minors will be to first place
a given apex graph's vertices on a wall, and then solve many routing
problems in order to show that the wall contains the appropriate
connections between these vertices. To do this, we need to partition
the one large wall into many smaller
walls. As shown in Figure~\ref{regions}, the union of the hexagons of a
wall forms a shape that can itself tile the plane, with a pattern of
connectivity equivalent to that of the original hexagonal tiling.
If the hexagons of a large wall are partitioned into smaller walls
according to such a tiling, we call the smaller walls {\em macrocells}.

Note that while the macrocells are a partition of the hexagons of a
wall, they are not a partition of the edges and vertices of the wall.
We say that two macrocells are {\em adjacent} if they share some edges
and vertices; if the macrocells are walls of size $s$ the shared
vertices form a path of $O(s)$ corner vertices (and possibly many more
path vertices). Define a {\em side} of a macrocell to be one of these
shared paths.

\begin{lemma}
\label{centralize}
Let $S$ be a set of $(s-t/2)$-inner corners of a wall of size $s$.
Then one can partition the wall into macrocells of size $t$ so that
$|S|/4$ members of $s$ are $t/2$-inner.
\end{lemma}

\begin{proof}
The partition into macrocells is determined by the choice of one central
hexagon for one macrocell.  If one chooses this hexagon uniformly at
random, the probability that any given corner in $S$ is $t/2$-inner is
propertional to the area of the inner size-$t/2$ wall of a macrocell
relative to the overall macrocell's area; this probability is therefore
$1/4$. Thus choosing a random macrocell center gives an expected number
of $t/2$-inner members of $S$ equal to $|S|/4$. The best macrocell
center must give at least as many $t/2$-inner members of $S$ as this
expectation.
\end{proof}

\begin{lemma}
\label{goodcurve}
Let $M$ be a set of macrocells of an embedded wall,
such that one can connect any two macrocells in $M$ by a chain of
adjacent pairs of macrocells.  Then there exists a non-self-intersecting
curve in the plane that is contained in the union of $M$, that passes
through all macrocells in
$M$, such that the intersection of the curve with any macrocell has at
most three connected components.
\end{lemma}

\begin{proof}
Form a planar graph by placing a point at the center of each macrocell,
and connecting pairs of points at the centers of adjacent macrocells.
Then by assumption this graph is connected, so we can choose a spanning
tree. A curve $C$ formed by thickening the edges of this tree and passing
around the boundary of the thickened tree has two of the three
properties we want: it is contained in the union of $M$ and passes
through all macrocells in $M$ (Figure~\ref{wallcurve}).

Now, suppose that some path $x$ in $C$
passes consecutively through three pairwise adjacent
macrocells $m_i$, $m_j$, and $m_k$, (e.g. at points where the spanning
tree edges form a
$60\dg$ angle), and the intersection of $C$ with the middle macrocell
$m_j$ has more than one component.
Then we can simplify $C$ by replacing $x$ with a curve that passes
directly from $m_i$ to $m_k$ (Figure~\ref{curvesimp}). This simplification
step maintains the two properties that $C$ is in the union of $M$ and
passes through each macrocell. It is possible for such a simplification
step to introduce a crossing, but only in the case that more than one
path passes through the same triple of macrocells; to avoid this problem
we always choose the innermost path when more than one path passes
through the same triple of macrocells. Each simplification step reduces
the total number of connected components formed by intersecting $C$ with
macrocells, so the simplification process must terminate.

Once this simplification process has terminated, the components of an
intersection of $C$ with a macrocell (if that intersection has multiple
components) must connect non-adjacent pairs of macrocells, so there can
be at most three components per macrocell.
\end{proof}

\section{Monotone Embedding}

We now show how to partition a graph into smaller pieces that can be
mapped onto a wall using $(s,t)$ routing problems.
Specifically, we will be concerned with performing this sort of
partition to walls, since any other planar graph can be found as a minor
of a wall (Lemma~\ref{wall-minor}).

\begin{defn}
Planar graph $G$ is {\em monotone embedded} in the plane if no vertical
line crosses any edge more than once, and no vertical line contains more
than one vertex.  The {\em monotone bandwidth} of $G$ is the maximum
number of edges crossed by any vertical line, minimized over all such
embeddings.
\end{defn}

\begin{lemma}\label{wall-band}
A wall of size $s$ has monotone bandwidth $O(s)$.
\end{lemma}

\begin{proof}
Draw the wall using regular hexagons, tilted slightly so that no edge is
vertical; this gives a monotone embedding. Any vertical line crosses
$O(s)$ hexagons, and hence $O(s)$ edges.
\end{proof}

\begin{figure}[p]
$$\efig{5in}{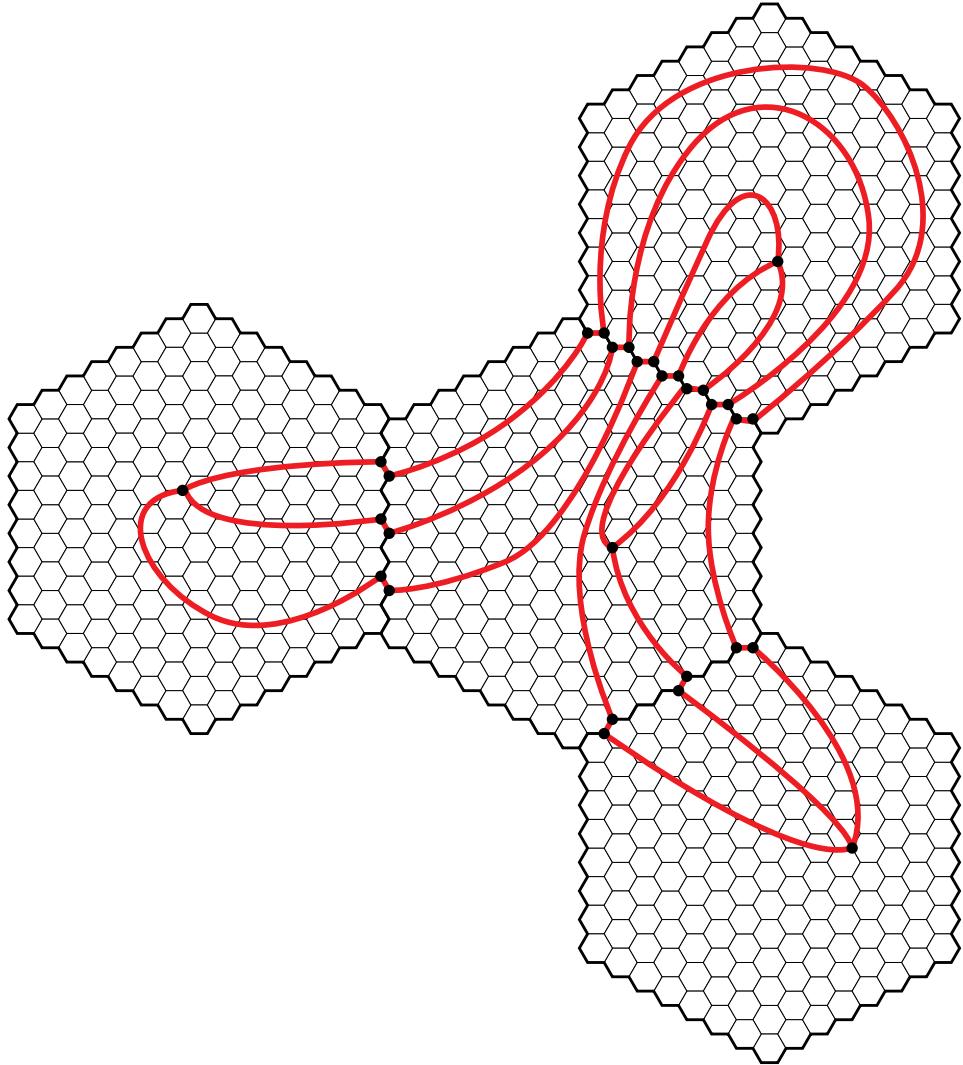}$$
\caption{The solutions to four routing problems on marked macrocells
can be combined to form a $K_4$ minor.}
\label{k4}
\end{figure}

\begin{lemma}\label{band-wall}
Let $W$ be a graph formed by connecting a sequence of macrocells of size
$s$, such that one can connect any two macrocells in $M$ by a chain of
adjacent pairs of macrocells, and let $k$ of the macrocells in $W$ have a
marked degree-three $s/2$-inner corner. Then $W$ contains as a
minor any $k$-vertex trivalent graph $G$ with monotone bandwidth
$s/6$, such that each subset of vertices of $W$ that is collapsed to
form each vertex of $G$ contains one of the marked vertices.
\end{lemma}

\begin{proof}
According to Lemma~\ref{goodcurve}, we can find a curve $C$ contained in
an embedding of $W$, and passing through each macrocell between one and
three times.
Find a monotone embedding of $G$, and form a correspondence with the
marked corners of $W$ (ordered by the positions along $C$ where $C$
first intersects each macrocell) and the vertices of $G$ (ordered
according to the monotone embedding).

Then we form an $(2,(s-2)/3)$ routing problem for each component of an
intersection of $C$ with a macrocell.  If the macrocell does not contain
a marked vertex, or if the component is not the first intersection of
$C$ with the macrocell, the routing problem just consists of pairs
of boundary vertices of the macrocell, with each pair placed on the two
sides of the macrocell crossed by $C$; the number of pairs is chosen to
match the number of edges cut by a vertical slice through the
corresponding part of the monotone embedding.  However, for the first
intersection of $C$ with a marked macrocell, we instead form a routing
problem in which the pattern of connections between the boundary vertices
and the marked inner corner matches the pattern of connections in a
vertical slice through the corresponding vertex of the monotone embedding.

This set of routing problems involves the placement of at most
$3+s/3$ terminals on any side of any macrocell.  These vertices can
be placed arbitrarily on that side, as long as they can be connected by
disjoint paths along the side to the corresponding terminals of the
adjacent macrocell.

The union of the at most three routing problems within each macrocell is
an $(s/2,s/2)$ routing problem and therefore has a solution.
Combining these solutions, and contracting the solution paths in each
macrocell, forms the desired minor.
\end{proof}

Figure~\ref{k4} depicts a set of routing problems on four macrocells,
the solutions to which could be combined to form a complete graph on
four vertices.  For simplicity we have drawn the figure using macrocells
of size eight, but (since $K_4$ has monotone bandwidth four) the lemma
above only guarantees the existence of such a routing for macrocells of
size 24.

We note that
Lemma~\ref{wall-minor} follows as an easy consequence of 
Lemma~\ref{band-wall}: given any $n$-vertex planar graph $G$, expand its
vertices into trees of degree-three vertices. The resulting $O(n)$-node
graph has monotone bandwidth $O(n)$, 
so it can be found as a minor of a
wall of size $O(n^{3/2})$, partitioned into a path of $O(n)$ smaller walls
of size $O(n)$ as depicted in Figure~\ref{regions}.

\section{The Main Result}

\begin{theorem}
\label{main}
Let $\cal F$ be a minor-closed family of graphs.
Then $\cal F$ has the diameter-treewidth property
iff $\cal F$ does not contain all apex graphs.
\end{theorem}

\begin{proof}
One direction is easy: we have seen that the apex graphs do not
have the diameter-treewidth property, so no family containing all apex
graphs can have the property.

In the other direction, we wish to show that if $\cal F$ does not have
the diameter-treewidth property, then it contains all apex graphs.  By
Lemma~\ref{wall-minor} it will suffice to find a graph in $\cal F$ formed by
connecting some vertex $v$ to all the vertices of a wall of size $n$,
for any given $n$.  If $\cal F$ does not have the diameter-treewidth
property, there is some $D$ such that $\cal F$ contains graphs with
diameter $D$ and with arbitrarily large treewidth.

Let $G$ be a graph in $\cal F$ with diameter $D$ and treewidth $W(N_1)$
for some large $N_1$ and for the function $W(N)$ shown to exist in
Lemma~\ref{wall}.  Then $G$ contains a wall of size $N_1$.  We choose
appropriate values $N_2$ and $N_3=\Theta(N_1/N_2)$ and partition the wall
into $N_3^2$ macrocells of size $N_2$. Say a macrocell is {\em good} if it
is not adjacent to the boundary of the wall.

Choose any vertex $v\in G$ and find a tree of shortest paths from
$v$ to each vertex.  We say that a macrocell is reached at level $i$ of
the tree if some vertex of the macrocell is included in that level.
Since $G$ has diameter $D$,the tree will have
height $D$. Since all macrocells are reached level $D$,
and the number of macrocells reached at level zero is just one,
there must be some intermediate level~$\lambda$ of the tree for which
the number $N_4$ of good macrocells reached is larger by a factor of
$N_3^{2/D}$ than the number of good macrocells reached in all
previous tree levels combined.

Let set $S$ be a set of corners of the wall formed by taking, in each
good macrocell reached at level $\lambda$, a corner nearest to one of the
vertices in that level of the tree.  By Lemma~\ref{centralize},
we can find a new partition into macrocells, and a set of $|S|/4$ corners
that are $N_2/2$-central for this partition. Each macrocell in this new
partition contains
$O(1)$ of these corners, so by removing corners that appear in the same
macrocell we can mark a set $S'$ of $\Omega(N_4)$
inner corners of macrocells, at most one corner per macrocell.
Note that the number of new macrocells reached at level $\lambda-1$
is still $O(N_4/N_3^{2/D})$, since each old macrocell reached at that
level can only contribute vertices to $O(1)$ new macrocells.

We then contract levels $1$ through $\lambda-1$ of the tree to a single
vertex $v$.  This gives a minor $G'$ of $G$ in which $v$ is connected
to inner corners of $\Omega(N_4)$ distinct macrocells, and in which
$O(N_4/N_3^{2/D})$ other macrocells are ``damaged'' by having a vertex
included in the contracted portion of the tree.  The adjacencies between
damaged regions of the wall form a planar graph with $O(N_4/N_3^{2/D})$
vertices and so $O(N_4/N_3^{2/D})$ faces, and there must therefore be a
face of this graph containing $\Omega(N_3^{2/D})$ members of $S'$. Let
$S''$ denote this subset of $S'$.

Now $S''$ is part of a connected set of undamaged macrocells of size
$N_2$, so by Lemma~\ref{band-wall} we can find a wall of size
$O(\min(N_3^{1/D},N_2)$ as a minor of this set of undamaged macrocells.
If $N_2=\Omega(n)$ and $N_3^{2/D}=\Omega(n^2)$, we can find a wall of
size $n$.  These conditions can both be assured by letting
$N_1=\Omega(n)^{D+1}$.  Combining this wall with the contracted vertex $v$
forms the apex graph minor we were seeking.

We can carry out this construction for any $n$,
and since by Lemma~\ref{wall-minor} every apex graph can be found as a
minor of graphs of the form of $M$, all apex graphs are minors of graphs
in $\cal F$ and are therefore themselves graphs of~$\cal F$.
\end{proof}

Alternately, instead of finding apex-grid graph minors, and using those
to find all other apex graphs as minors, we can find any apex graph
directly by following the proof of Lemma~\ref{wall-minor} sketched above
after Lemma~\ref{band-wall}.

\section{Bounded Genus Graphs}

The results above show that any minor-closed family excluding an apex
graph has the diameter-treewidth property.  For example, consider the
bounded genus graphs. It is not hard to show that, for any
$g$, there is an apex graph with genus more than $g$: genus $g$ graphs
have at most $3n+O(g)$ edges, while maximal apex graphs have $4n-10$
edges, so choosing $n$ large gives an apex graph with too many edges to
have genus~$g$. Therefore, genus $g$ graphs have the diameter-treewidth
property. However this proof does not give us a very tight relation
between diameter, genus, and treewidth. We can achieve a much better
treewidth bound by proving the diameter-treewidth property more directly.

\begin{figure}
$$\efig{4in}{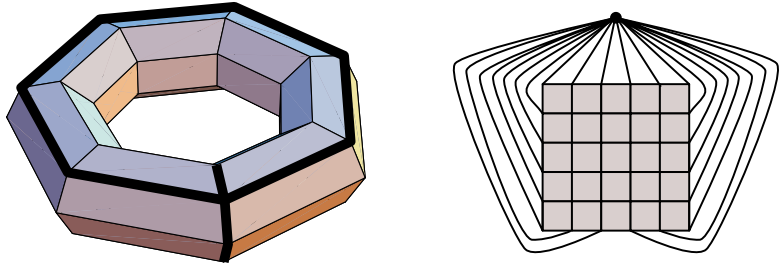}$$
\caption{Torus graph with subgraph $X$ highlighted, and planar graph
formed by contracting $X$}
\label{genus}
\end{figure}

\begin{lemma}
\label{splitter}
Let $G$ be embedded on a surface $S$ of genus $g$, with all faces
of the embedding topologically equivalent to disks. Then there exists a
subgraph
$X$ of $G$, isomorphic to a subdivision of a graph $Y$ with $O(g)$ edges
and vertices, such that the removal of the
points of $X$ from $S$ leaves a set  topologically equivalent to a disk.
\end{lemma}

\begin{proof}
Let $X$ be a minimal connected subgraph of $G$ such that all components of
$S-X$ are topological disks.  Then there must be at most one such
component, for multiple components could be merged by
removing from $X$ an edge along which two adjacent components are
connected; any such merger preserves the disk topology of the components
and the connectivity of $X$ (since any path through the removed edge can
be replaced by a path around the boundary of a component).

Thus $X$ is a graph bounding a single disk face.
By Euler's formula, if $X$ has $n$ vertices, it has $n+O(g)$ edges.
Let $T$ be a spanning tree of $X$; then $X-T$ has $O(g)$ edges.
Note also that $X$ has no degree-one vertices,
so each leaf of $T$ must be an endpoint of an edge in $X-T$
and there are $O(g)$ leaves.
Any graph formed by adding $O(g)$ edges to a tree with $O(g)$ leaves
must be a subdivision of a graph with $O(g)$ edges and vertices.
\end{proof}

Figure~\ref{genus} depicts a graph $X$ for an example in which
$G$ is embedded on a torus.

\begin{theorem}
Let $G$ have genus $g$ and diameter $D$.
Then $G$ has treewidth $O(gD)$.
\end{theorem}

\begin{proof}
Embed $G$ on a minimal-genus surface $S$, so that all its faces are
topological disks. Choose a subgraph $X$ as in Lemma~\ref{splitter},
having the minimum number of edges possible among all subgraphs satisfying
the conditions of the lemma, and let $Y$ be a graph with $O(g)$
vertices and edges of which $X$ is a subdivision (as described in the
lemma). Then, each path in $X$ corresponding to an edge
in $Y$ has $O(D)$ edges. For, if not, one could find a smaller $X$ by
replacing part of a long path by the shortest path from its midpoint to
the rest of $X$. Therefore, $X$ has $O(gD)$ edges and vertices.

Now contract $X$ forming a minor $G'$ of $G$. The result is a planar
graph, since $G-X$ can remain embedded in its disk, with the
vertex contracted from $X$ being connected to $G-X$ by edges that cross
the boundary of this disk. The contraction can only reduce the diameter
of $G$. Therefore,
$G'$ has treewidth $O(D)$, and a tree decomposition of $G$ with
treewidth $O(gD)$ can be formed by adjoining $X$ to each clique in a
tree decomposition of $G'$.
\end{proof}

\section{Algorithmic Consequences}

\begin{theorem}
For any minor-closed family of graphs with the
diameter-treewidth property, there exists a linear time
approximation scheme for maximimum independent set, minimum vertex cover,
maximum $H$-matching, minimum dominating set, and the other approximation
problems solved by Baker~\cite{Bak-JACM-94}.
\end{theorem}

The method is the same as in \cite{Bak-JACM-94}: we remove every $k$th
level in a breadth first search tree, with one of $k$ different choices
of the starting level, forming a collection of subgraphs each of which
is induced by some $k-1$ contiguous levels of the tree. (Forthe minimum
dominating set and vertex cover problems, we instead duplicate the
vertices on every $k$th level, and form subgraphs induced by $k+1$
contiguous levels of the tree). As Baker shows, one of these choices leads
to a graph that approximates the optimum within a $1+O(1/k)$ factor.
We then use the diameter-treewidth property to show that each of these
subgraphs has bounded treewidth. A tree decomposition of each subgraph
can be found in linear time \cite{Bod-STOC-93}, after which the
appropriate optimization problem can be solved in linear time in each
subgraph by using dynamic programming techniques \cite{BerLawWon-Algs-87,
TakNisSai-JACM-82}.

We note that maximum independent set
can also be approximated for all minor-closed families, using the results
of Alon et al. \cite{AloSeyTho-STOC-90} on separator theorems for such
families, however the separator algorithm of \cite{AloSeyTho-STOC-90}
takes superlinear time $O(k^{1/2}n^{3/2})$ (where $k$ is the number of
vertices of the largest clique belonging to the family) and this
approximation technique does not seem to apply to the other problems on
the list above.

\begin{theorem}
Subgraph isomorphism or induced subgraph isomorphism for a fixed pattern
$H$ in any minor-closed family of graphs with the diameter-treewidth
property can be tested in time $O(n)$.
\end{theorem}

The algorithm closely follows that of \cite{Epp-JGAA-97}.
We again remove every $k$th level of the tree with one of $k$ different
choices of the starting level, forming subgraphs of $k-1$ contiguous
levels, where $k-1$ is the diameter of $H$. If $H$ occurs in $G$, it
must occur in one of these subgraphs, which can be tested by finding a
tree decompostion and performing dynamic programming.

\section{Conclusions and Open Problems}

We have characterized
the minor-closed families with the diameter-treewidth property.
However, some further work remains.
Notably, the relation we showed between diameter and treewidth
was not as strong as for planar graphs: for planar graphs (and
bounded-genus graphs)
$w=O(d)$ while for other minor-closed families our proof only shows that
$w=W(c^{d+1}))$, where $c$ is a constant that depends on the family and
$W(x)$ represents the rapidly-growing function used by Robertson and
Seymour to prove Lemma~\ref{wall}.  Can we prove tighter bounds on
treewidth for general minor-closed families?

Specifically, what relation between diameter and treewidth holds for
the graphs having no $K_{3,a}$ minor for some fixed $a$?  Note that
$K_{3,a}$ is an apex graph, so these graphs have the diameter-treewidth
property. $K_{3,a}$-free graphs are a generalization of
planar graphs (which have no $K_{3,3}$ or $K_5$ minor)
and have other interesting properties;
notably, in connection with the subgraph isomorphism algorithms
described above, a subgraph $H$ has an $O(n)$ bound on the number of
times it can occur in $K_{3,a}$-free graphs, if and only if $H$ is
3-connected~\cite{Epp-JGT-93}. Any improved treewidth bounds would
improve the running time and practicality of the subgraph isomorphism and
approximation algorithms we described.

Also, are there natural families of graphs that are not minor-closed and
that have the diameter-treewidth property (other than the bounded-degree
graphs or other classes in which a diameter bound imposes a limit on
total graph size)?  Although one could not then apply Baker's
approximation technique \cite{Bak-JACM-94}, this would still lead to
quadratic-time subgraph isomorphism algorithms based on testing
bounded-radius neighborhoods of each vertex \cite{Epp-JGAA-97}.

Finally, can we extend some of the same efficient subgraph isomorphism
and approximation algorithms to graph families without the
diameter-treewidth property? For instance, it is trivial to do so for
apex graphs, by treating the apex specially and applying a modified
algorithm in the remaining graph. What about other graph families
containing the apex graphs, such as linkless and knotless embeddable
graphs, or $K_{4,4}$-free graphs?

\frenchspacing
\bibliography{subiso}
\end{document}